\documentclass[a4paper,reqno]{amsart}%
\usepackage{amssymb}
\usepackage{amsmath}
\usepackage{amsfonts}
\usepackage{graphicx}
\usepackage{amscd}%
\theoremstyle{plain}

\newtheorem{lemma}{Lemma}

\numberwithin{equation}{section}

\begin{document}
\title[Erratum: Quantum Montgomery identity and Ostrowski type inequalities]{Erratum: Quantum Montgomery identity and quantum estimates of Ostrowski type
inequalities}
\author[A.~Aglić~Aljinović]{Andrea Agli\'{c} Aljinovi\'{c}}
\address{University of Zagreb, Faculty of Electrical Engineering and Computing, Unska
3, 10 000 Zagreb, Croatia}
\email{andrea.aglic@fer.hr}
\author[D.~Kovačević]{Domagoj Kova\v{c}evi\'{c}}
\address{University of Zagreb, Faculty of Electrical Engineering and Computing, Unska
3, 10 000 Zagreb, Croatia}
\email{domagoj.kovacevic@fer.hr}
\author[M.~Puljiz]{Mate Puljiz}
\address{University of Zagreb, Faculty of Electrical Engineering and Computing, Unska
3, 10 000 Zagreb, Croatia}
\email{mate.puljiz@fer.hr}
\date{July 17, 2020}

\begin{abstract}
We disprove and correct some recently obtained results regarding Montgomery
identity for quantum integral operator and Ostrowski type inequalities
involving convex functions.
\end{abstract}

\maketitle

\section{Introduction}

In \cite{KUNT} the authors obtained the following
generalization of Montgomery identity for quantum calculus.\smallskip

\begin{lemma}
\cite{KUNT} (Quantum Montgomery identity) Let $f:\left[  a,b\right]
\rightarrow%
\mathbb{R}
$, be an arbitrary function with $d_{q}^{a}f$ quantum integrable on $[a,b]$,
then the following quantum identity holds:
\begin{equation}
f\left(  x\right)  -\frac{1}{b-a}%
{\displaystyle\int\limits_{a}^{b}}
f\left(  t\right)  d_{q}^{a}t=\left(  b-a\right)
{\displaystyle\int\limits_{0}^{1}}
K_{q,x}\left(  t\right)  D_{q}^{a}f\left(  tb+\left(  1-t\right)  a\right)
d_{q}^{0}t\label{qMI}%
\end{equation}
where $K_{q,x}\left(  t\right)  $ is defined by
\begin{equation}
K_{q,x}\left(  t\right)  =\left\{
\begin{array}
[c]{cc}%
qt, & 0\leq t\leq\dfrac{x-a}{b-a},\\
& \\
qt-1, & \dfrac{x-a}{b-a}<t\leq1.
\end{array}
\right.  \label{PK}%
\end{equation}
\smallskip
\end{lemma}

Using this identity, the authors have obtained two Ostrowski type
inequalities for quantum integrals and applied it in several special cases.
\smallskip

Unfortunately, in the proof of this lemma an error is made when calculating the
integrals involving the kernel $K_{q,x}\left(  t\right)  $ on the interval
$\left[  \frac{x-a}{b-a},1\right]  $. Also, in the proofs of Theorem 3 and
Theorem 4 a small mistake related to the convexity of $\left\vert
D_{q}^{a}f\right\vert ^{r}$ is made. \medskip

In the present paper we prove that the identity (\ref{qMI}) and, thus, all
of the consequent results are incorrect and provide corrections for these
results.\bigskip

\section{Main results\medskip}

The q-derivative
of a function $f:\left[  a,b\right]  \rightarrow%
\mathbb{R}
$ for $q\in\left\langle 0,1\right\rangle $
(see \cite{TARI} or \cite{JACK} for $a=0$) is given by
\begin{align*}
D_{q}^{a}f\left(  x\right)   &  =\frac{f\left(  x\right)  -f\left(  a+q\left(
		x-a\right)  \right)  }{\left(  1-q\right)  \left(  x-a\right)  },\text{ for
}x\in\left\langle a,b\right]  \text{\ }\bigskip\\
D_{q}^{a}f\left(  a\right)   &  =\underset{x\rightarrow a}{\lim}\text{ }%
D_{q}^{a}f\left(  x\right)  \text{\ }\bigskip
\end{align*}
We say that $f:\left[  a,b\right]  \rightarrow%
\mathbb{R}
$ is q-differentiable
if $\underset{x\rightarrow a}{\lim}D_{q}^{a}f\left(  x\right)  $ exists. The q-derivative is a discretization of the ordinary
derivative and if $f$ is a differentiable function then (\cite{ANNA},
\cite{KAC})
\[
\underset{q\rightarrow1}{\lim}\text{ }D_{q}^{a}f\left(  x\right)  =f^{\prime
}\left(  x\right)  .
\]
Further, the q-integral of $f$ is defined by
\[%
{\displaystyle\int\limits_{a}^{x}}
f\left(  t\right)  d_{q}^{a}t=\left(  1-q\right)  \left(  x-a\right)
{\displaystyle\sum\limits_{k=0}^{\infty}}
q^{k}f\left(  a+q^{k}\left(  x-a\right)  \right)  ,\text{ \ }x\in\left[
a,b\right]  .
\]
If the series on the right hand-side is convergent, then the q-integral $%
{\displaystyle\int_{a}^{x}}
f\left(  t\right)  d_{q}^{a}t$ exists and $f:\left[  a,b\right]  \rightarrow%
\mathbb{R}
$ is said to be q-integrable on $\left[  a,x\right]  $.
If $f$ is continuous on $\left[
a,b\right]  $ the series $\left(  1-q\right)  \left(
x-a\right)
{\displaystyle\sum\limits_{k=0}^{\infty}}
q^{k}f\left(  a+q^{k}\left(  x-a\right)  \right)  $ tends to the Riemann
integral of $f$ as $q\rightarrow1$ (\cite{ANNA}, \cite{KAC})
\[
\underset{q\rightarrow1}{\lim}%
{\displaystyle\int\limits_{a}^{x}}
f\left(  t\right)  d_{q}^{a}t=%
{\displaystyle\int\limits_{a}^{x}}
f\left(  t\right)  dt.
\]
If $c\in\left\langle a,x\right\rangle $ the q-integral is defined by
\[%
{\displaystyle\int\limits_{c}^{x}}
f\left(  t\right)  d_{q}^{a}t=%
{\displaystyle\int\limits_{a}^{x}}
f\left(  t\right)  d_{q}^{a}t-%
{\displaystyle\int\limits_{a}^{c}}
f\left(  t\right)  d_{q}^{a}t.
\]
Obviously, the q-integral depends on the values of $f$ at the points outside the
interval of integration and an important difference between the definite
q-integral and Riemann integral is that even if we are integrating a function
over the interval $[c,x]$, $a<c<x<b$, for q-integral
we have to take into account its behavior at
$t=a$ as well as its values on $[a,x]$. This is the main reason for mistakes
made in \cite{KUNT} since in the proof of Lemma 1 the following error was made:
\begin{align*}
&
{\displaystyle\int\limits_{\frac{x-a}{b-a}}^{1}}
K_{q,x}\left(  t\right)  D_{q}^{a}f\left(  tb+\left(  1-t\right)  a\right)
d_{q}^{0}t\\
&  =%
{\displaystyle\int\limits_{0}^{1}}
\left(  qt-1\right)  D_{q}^{a}f\left(  tb+\left(  1-t\right)  a\right)
d_{q}^{0}t-%
{\displaystyle\int\limits_{0}^{\frac{x-a}{b-a}}}
\left(  qt-1\right)  D_{q}^{a}f\left(  tb+\left(  1-t\right)  a\right)
	d_{q}^{0}t.
\end{align*}
But $K_{q,x}\left(  t\right)  \neq\left(  qt-1\right)  $ for $t\in\left[
0,1\right]  $ or for $t\in\left[  0,\frac{x-a}{b-a}\right]  $,
so the equality does not hold.\medskip

Now, we give a proof that the quantum Montgomery identity (\ref{qMI}) is not
correct, since it does not hold for all $x\in\left[  a,b\right]  $. As we
shall see, the identity (\ref{qMI}) is valid only if $x=a+q^{m+1}\left(
b-a\right)  $ for some $m\in%
\mathbb{N}
\cup\left\{  0\right\}  $. We have%
\begin{align*}
&  \left(  b-a\right)
{\displaystyle\int\limits_{0}^{1}}
K_{q,x}\left(  t\right)  D_{q}^{a}f\left(  tb+\left(  1-t\right)  a\right)
d_{q}^{0}t\\
&  =\left(  b-a\right)  \left(  1-q\right)
{\displaystyle\sum\limits_{k=0}^{\infty}}
q^{k}K_{q,x}\left(  q^{k}\right)  D_{q}^{a}f\left(  a+q^{k}\left(  b-a\right)
\right)  .
\end{align*}
For $q\in\left\langle 0,1\right\rangle $ let $m\in%
\mathbb{N}
\cup\left\{  0\right\}  $ be such that
\[
q^{m+1}\leq\frac{x-a}{b-a}<q^{m},
\]
in other words
\[
m=\left\lceil \log_{q}\frac{x-a}{b-a}\right\rceil -1.
\]
Then
\[
K_{q,x}\left(  q^{k}\right)  =\left\{
\begin{array}
[c]{cl}%
q^{k+1}-1, & k\leq m,\\
& \\
q^{k+1}, & k\geq m+1,
\end{array}
\right.
\]
and%
\begin{align*}
&  \left(  b-a\right)  \left(  1-q\right)
{\displaystyle\sum\limits_{k=0}^{\infty}}
q^{k}K_{q,x}\left(  q^{k}\right)  D_{q}^{a}f\left(  a+q^{k}\left(  b-a\right)
\right)  \\
&  =\left(  b-a\right)  \left(  1-q\right)  \left(
{\displaystyle\sum\limits_{k=0}^{m}}
q^{k}\left(  q^{k+1}-1\right)  \frac{f\left(  a+q^{k}\left(  b-a\right)
\right)  -f\left(  a+q^{k+1}\left(  b-a\right)  \right)  }{\left(  1-q\right)
q^{k}\left(  b-a\right)  }\right.  \\
&  \left.  +%
{\displaystyle\sum\limits_{k=m+1}^{\infty}}
q^{k}\left(  q^{k+1}\right)  \frac{f\left(  a+q^{k}\left(  b-a\right)
\right)  -f\left(  a+q^{k+1}\left(  b-a\right)  \right)  }{\left(  1-q\right)
q^{k}\left(  b-a\right)  }\right)  \\
&  =-%
{\displaystyle\sum\limits_{k=0}^{m}}
\left(  f\left(  a+q^{k}\left(  b-a\right)  \right)  -f\left(  a+q^{k+1}%
\left(  b-a\right)  \right)  \right)  \\
&  +%
{\displaystyle\sum\limits_{k=0}^{\infty}}
q^{k+1}\left(  f\left(  a+q^{k}\left(  b-a\right)  \right)  -f\left(
a+q^{k+1}\left(  b-a\right)  \right)  \right)  \\
&  =f\left(  a+q^{m+1}\left(  b-a\right)  \right)  -f\left(  b\right)  +%
{\displaystyle\sum\limits_{k=0}^{\infty}}
q^{k+1}\left(  f\left(  a+q^{k}\left(  b-a\right)  \right)  -f\left(
a+q^{k+1}\left(  b-a\right)  \right)  \right)  .
\end{align*}
If we put $S=%
{\displaystyle\sum\limits_{k=0}^{\infty}}
q^{k}f\left(  a+q^{k}\left(  b-a\right)  \right)  =\frac{1}{\left(
1-q\right)  \left(  b-a\right)  }%
{\displaystyle\int\limits_{a}^{b}}
f\left(  t\right)  d_{q}^{a}t$, we have%
\[%
{\displaystyle\sum\limits_{k=0}^{\infty}}
\left(  q^{k+1}\right)  \left(  f\left(  a+q^{k}\left(  b-a\right)  \right)
-f\left(  a+q^{k+1}\left(  b-a\right)  \right)  \right)  =qS-\left(
S-f\left(  b\right)  \right)
\]
and
\begin{align*}
&  \frac{1}{b-a}%
{\displaystyle\int\limits_{a}^{b}}
f\left(  t\right)  d_{q}^{a}t+\left(  b-a\right)
{\displaystyle\int\limits_{0}^{1}}
K_{q,x}\left(  t\right)  D_{q}^{a}f\left(  tb+\left(  1-t\right)  a\right)
d_{q}^{0}t\\
&  =\frac{1}{b-a}%
{\displaystyle\int\limits_{a}^{b}}
f\left(  t\right)  d_{q}^{a}t+\left(  f\left(  a+q^{m+1}\left(  b-a\right)
\right)  -f\left(  b\right)  \right)  +qS-\left(  S-f\left(  b\right)
\right)  \\
&  =\left(  1-q\right)  S+f\left(  a+q^{m+1}\left(  b-a\right)  \right)
-f\left(  b\right)  +qS-S+f\left(  b\right)  \bigskip\\
&  =f\left(  a+q^{m+1}\left(  b-a\right)  \right)
\end{align*}
which is obviously not equal to $f\left(  x\right)  $, unless $x=a+q^{m+1}%
\left(  b-a\right)  $.\bigskip

This is no surprise since
Jackson integral takes into account only $f\left(  a+q^{k}\left(  x-a\right)
\right)  $ for $k\in%
\mathbb{N}
\cup\left\{  0\right\}  $. Thus, we have proved the next
lemma which is a corrected version of Lemma 1 from \cite{KUNT}\textbf{.}

\begin{lemma}
\label{1} (Quantum Montgomery identity) Let $f:\left[  a,b\right]  \rightarrow%
\mathbb{R}
$, be an arbitrary function with $D_{q}^{a}f$ quantum integrable on $[a,b]$,
then for all $x\in\left\langle a,b\right\rangle $ the following quantum
identity holds:
\begin{align*}
&  f\left(  a+q^{\left\lceil \log_{q}\frac{x-a}{b-a}\right\rceil }\left(
b-a\right)  \right)  -\frac{1}{b-a}%
{\displaystyle\int\limits_{a}^{b}}
f\left(  t\right)  d_{q}^{a}t\\
&  =\left(  b-a\right)
{\displaystyle\int\limits_{0}^{1}}
K_{q,x}\left(  t\right)  D_{q}^{a}f\left(  tb+\left(  1-t\right)  a\right)
d_{q}^{0}t
\end{align*}
where $K_{q,x}\left(  t\right)  $ is defined by
\[
K_{q,x}\left(  t\right)  =\left\{
\begin{array}
[c]{cc}%
qt, & 0\leq t\leq\dfrac{x-a}{b-a},\\
& \\
qt-1, & \dfrac{x-a}{b-a}<t\leq1.
\end{array}
\right.
\]

\end{lemma}

In Theorem 3 and Theorem 4 from \cite{KUNT} the authors have used the identity
(\ref{qMI}) to derive Ostrowski type inequalities for functions $f$ for which
$D_{q}^{a}f$ \ is quantum integrable on $[a,b]$ and $\left\vert D_{q}%
^{a}f\right\vert ^{r}$, $r\geq1$ is a convex function. Since these inequalities
depends on the validity of Lemma 1, our discussion invalidates all the results
from \cite{KUNT}.\smallskip

More precisely, in all the inequalities \linebreak $f\left(  x\right)$ should be
substituted with
$f\left(  a+q^{\left\lceil \log_{q}\frac{x-a}{b-a}\right\rceil }\left(
b-a\right)  \right)  $. In Theorems 3 and 4 $\left\vert D_{q}^{a}f\left(
a\right)  \right\vert ^{r}$ and $\left\vert D_{q}^{a}f\left(  b\right)
\right\vert ^{r}$ should be swapped, since in the proofs of Theorem 3 and Theorem
4, when applying the convexity of $\left\vert D_{q}^{a}f\right\vert ^{r}$
the following mistake was made
\[
\left\vert D_{q}^{a}f\left(  tb+\left(  1-t\right)  a\right)  \right\vert
^{r}\leq t\left\vert D_{q}^{a}f\left(  a\right)  \right\vert ^{r}+\left(
1-t\right)  \left\vert D_{q}^{a}f\left(  b\right)  \right\vert ^{r}.
\]
\medskip

\section*{Acknowledgments}
Domagoj Kova\v{c}evi\'{c} was supported by the QuantiXLie Centre of Excellence, a project
co financed by the Croatian Government and European Union through the
European Regional Development Fund - the Competitiveness and Cohesion
Operational Programme (Grant KK.01.1.1.01.0004).

\end{document}